\documentclass[12pt,fleqn]{article}
\usepackage{amsmath,amssymb,amsthm,enumerate}%,backref,cite
\usepackage[cp1251]{inputenc}
\usepackage[ukrainian]{babel}
\usepackage[T2A]{fontenc}
\usepackage[all]{xy}

\begin{document}
УДК 517.9
\bigskip

\textbf{Ю. Б. Зелинский, Б. А. Клищук, М. В. Ткачук} (Yu. B.
Zelinskii, B. A. Klishchuk, M. V. Tkachuk), Ін-т математики НАН
України, Київ.
\bigskip

\textbf{ТЕОРЕМЫ О ВКЛЮЧЕНИИ ДЛЯ МНОГОЗНАЧНЫХ ОТОБРАЖЕНИЙ}

(\textbf{THEOREMS ABOUT INCLUDINGS FOR MULTIVALUED MAPPINGS}).
\bigskip

This paper is devoted to studying of some properties of multivalued
mappings in Euclidean space. There were proved theorems on a fixed
point for multi\-valued mappings whose restrictions to some subset
in the closure of a domain satisfy ``a coacute angle condition''\/
or ``a strict coacute angle condition''\/. There also were obtained
similar results for restrictions of multivalued mappings satisfying
some metric limitations.
\medskip

Вивчаються деякі властивості многозначних відображень в евклідовому
просторі. Доведено теореми про нерухому точку для многозначних
відображень, звуження яких на деяку підмножину в замиканні області
задовольняють ``умові когострого кута''\/ або ``умові строгого
когострого кута''\/. Подібні результати отримано і для звужень
многозначних відображень, які задовольняють деяким метричним
обмеженням.
\medskip

Изучаются некоторые свойства многозначных отображений в евклидовом
пространстве. Доказаны теореми о неподвижной точке для многозначных
отображений, сужения которых на некоторое подмножество в замыкании
области удовлетворяют ``условию коострого угла''\/ либо ``условию
строгого коострого угла''\/.  Подобные результаты получены и для
сужений многозначных отображений, удовлетворяющих некоторым
метрическим ограничениям.

\newpage

{\bf 1. Введение.} В настоящей работе мы продолжаем исследование
многозначных включений, начатых в [1], и основанных на использовании
геометрической формы теоремы Хана-Банаха. Мы избавляемся от условия
содержания искомой областью начала координат, рассматриваем
отображения не только в то же пространство, но и в другое, а также
уменьшаем размеры множества, на котором справедливы ``ограничения
типа острого угла'' [2-5].

{\bf 2. Обозначения и основные определения.} Пусть  $E^{n}$ ---
$n$-мерное евклидово (действительное или комплексное) пространство,
$x$, $y$ --- некоторые точки $ E ^ {n} $, $A$, $B$ --- подмножества
$ E ^ {n} $, $\langle\ast ,\ast\rangle$
--- скалярное произведение в $E^{n}$, $conv\ A$ --- выпуклая
оболочка множества $A$, $\angle xOy = arccos\frac{Re\langle x,
y\rangle}{\sqrt{\langle x, x\rangle}\sqrt{\langle y, y\rangle}}$.

Далее будем рассматривать многозначные (в том числе одно\-знач\-ные
и разрывные) отображения подмножеств евклидового про\-странства.

Пусть $X$ и $Y$ --- топологические пространства. Если $F_{1}, F_{2}:
X \rightarrow Y$ --- два многозначных отображения, то будем
говорить, что $F_{2}$ есть сужением отображения $F_{1}$ , если
$F_{1}(x)\supset F_{2}(x)\neq \varnothing$ для всех точек $x\in X$(в
частности, если $A\supset B$ и $F_{1}:A\rightarrow Y$,
$F_{2}:B\rightarrow Y$ --- два отображения, то отображение $F_{2}$
есть сужением отображения $F_{1}$ на $B$, если $F_{1}(x) = F_{2}(x)$
при $x \in B$  и $F_{2}(x) = \varnothing$ при $x \notin B$, т.е. не
исключено, что для отдельных точек образы сужения пустые множества).

Скажем, что на множестве $A$ отображение $F$ удовлетворяет ``условию
острого (строгого острого) угла'', если выполнено условие ${\rm
Re}\, \langle x, y \rangle\geq 0$\ (${\rm Re}\, \langle x, y
\rangle>0$) для всех пар точек $x\in A$, $y\in F(x)$.

{\bf 3. ``Условие коострого угла''.} Пусть $Y^{\ast}$  --- дуальное
пространство к пространству $Y$. Будем говорить, что отображение $F:
A\rightarrow Y$ $(A\subset X)$ удовлетворяет "условию коострого
угла" \ на $A$, если для каждой точки $y^{\ast}\in Y^{\ast}$,
$y^{\ast}\neq 0$, существует точка $x\in A$ такая, что выполнено
условие ${\rm Re}\, \langle y, y^{\ast} \rangle\geq 0$ для всех
точек $y\in F(x)$.

%\newpage

Справедливо следующее утверждение.

{\bf Теорема 1.} {\it Пусть $D$
--- область  евклидова пространства $E^{n} = X$. Пусть $K \subset
\overline{D}$ --- под\-мно\-жество в замыкании этой области и пусть
существует такое сужение $F_{1}$ многозначного отображения $F:
\overline{D}\rightarrow E^{n} = Y$ на подмножество $K$, которое
удовлетворяет ``условию коострого угла''  и $conv\ F_{1}(K)$ ---
компакт. Тогда если $conv\ F_{1}(K) \subset F(\overline{D})$, то $0
\in F(\overline{D})$.}

\textbf{\textit {Доказательство.}} Предположим, что $0 \notin
F(\overline{D})$. Следовательно, $0 \notin conv\ F_{1}(K)$. Тогда
согласно геометрической форме теоремы Хана--Банаха существует
гиперплоскость $L$, которая отделяет начало координат от компактного
выпуклого множества $conv\ F_{1}(K)$. Выберем луч $l$, выходящий из
начала координат перпендикулярно к гиперплоскости $L$ и направленный
в сторону противоположную $conv\ F_{1}(K)$. Для евклидовых
пространств отображение двойственности $\mathfrak{F}:Y \rightarrow
Y^{\ast}$ биективно. Выберем произвольную точку $y^{\ast}$, отличную
от начала координат, на луче $l$. С одной стороны по построению
$y^{\ast}\notin conv\ F_{1}(K)$, а с другой, согласно ``условию
коострого угла'', существует точка $x\in K$, образ $F_{1}(x)$
которой должен находиться в том же полупространстве по отношению к
гиперплоскости $L$, что и точка $y^{\ast}$. Полученное противоречие
доказывает теорему.

{\bf 3. ``Условие  строгого коострого угла''.} Отображение $F:
A\rightarrow Y$ $(A\subset X)$ удовлетворяет "условию строгого
коострого угла" \ на $A$, если для каждой точки $y^{\ast}\in
Y^{\ast}$, $y^{\ast}\neq 0$, существует точка $x\in A$ такая, что
выполнено условие ${\rm Re}\, \langle y,y^{\ast} \rangle>0$ для всех
точек $y\in F(x)$.

Используя рассуждения предыдущей теоремы вместе с рассуждениями,
примененными при доказательстве теоремы 2 [1], получим следующий
результат.

\newpage

{\bf Теорема 2.} {\it Пусть $D$  --- область  евклидова пространства
$E^{n} = X$. Пусть $K \subset \overline{D}$
--- под\-мно\-жество в замыкании этой области и пусть существует такое
сужение $F_{1}$ многозначного отображения $F:
\overline{D}\rightarrow E^{n} = Y$ на подмножество $K$, которое
удовлетворяет ``условию строгого коострого угла''. Тогда если $conv\
F_{1}(K) \subset F(\overline{D})$, то $0 \in F(\overline{D})$.}

\textbf{\textit {Доказательство.}} Предположим, что $0 \notin
F(\overline{D})$ и, следовательно, $0 \notin conv\ F_{1}(K)$.
Внутренность $Int\ (conv\ F_{1}(K))$ будет выпуклым открытым
множеством, не содержащим начало координат. Если $Int\ (conv\
F_{1}(K))= \varnothing$, то множество $conv\ F_{1}(K)$ имеет
размерность не выше $(n-1)$  и поэтому полностью лежит в некоторой
гиперплоскости. Следовательно существует гиперплоскость $L$, которая
проходит через начало координат и которая или полностью содержит
множество $conv\ F_{1}(K)$, или же с ним не пересекается. Если же
$Int\ (conv\ F_{1}(K))\neq \varnothing$, то существует
гиперплоскость $L$, которая проходит через начало координат и не
пересекает множество $Int\ (conv\ F_{1}(K))$. Для произвольного
выпуклого множества $A$ с непустой внутренностью ($Int\ A\neq
\varnothing$) справедливо $\overline{Int\ A}= \overline{A}$.
Следовательно, в обоих случаях множество $conv\ F_{1}(K)$ полностью
лежит в одном из замкнутых полупространств, на которые плоскость $L$
разбивает пространство. Теперь можем выбрать луч $l$, выходящий из
начала координат перпендикулярно к гиперплоскости $L$ и направленный
в сторону, противоположную полупространству, содержащему множество
$conv\ F_{1}(K)$. Выберем произвольную точку $y^{\ast}\in l$,
отличную от начала координат, на этом луче. С одной стороны
$y^{\ast}\notin conv\ F_{1}(K)$ , а с другой, согласно ``условию
строгого коострого угла'', существует точка $x\in K$ , образ
$F_{1}(x)$ которой должен находиться в том же полупространстве по
отношению к гиперплоскости $L$, что и точка $y^{\ast}$. Полученное
противоречие доказывает теорему.

{\bf 4. ``$\varepsilon$-условие острого угла''.} Дальше зададимся
целью уменьшить размеры множества, на котором выполнены ``условия
типа острого угла'' за счет более строгих неравенств.

Скажем, что множество $A$ является радианной (угловой)
$\varepsilon$-сетью, если для каждого луча, выходящего из начала
координат, существует луч, образующий с ним угол радианной величины
меньше $\varepsilon$  и пересекающий $A$.

Будем говорить, что на множестве $A$ отображение $F$ удовлетворяет
``$\varepsilon$-условию острого угла'', если $X = Y$ и для
произвольной точки $x\in A$ существует $y\in F(x)$ такое, что
выполняется условие $\angle xOy < \frac{\pi}{2}-\varepsilon$.

{\bf Теорема 3.} {\it Пусть $D$ --- область  евклидова пространства
$E^{n} = X$. Пусть $K \subset \overline{D}$ --- подмножество в
замыкании этой области, являющееся радианной $\varepsilon$-сетью и
пусть существует сужение $F_{1}$ многозначного отображения $F:
\overline{D}\rightarrow E^{n} = X$  на подмножество $K$, которое
удовлетворяет ``$\varepsilon$-условию острого угла''. Тогда, если
$conv\ F_{1}(K) \subset F(\overline{D})$, то $0 \in
F(\overline{D})$.}

\textbf{\textit {Доказательство.}} Предположим, что $0 \notin
F(\overline{D})$ . Как и в теореме 2 найдем гиперплоскость $L$,
которая проходит через начало координат и которая или полностью
содержит множество $conv\ F_{1}(K)$ , или же с ним не пересекается.
Выберем луч $l$, выходящий из начала координат перпендикулярно к
гиперплоскости $L$  и направленный в сторону противоположную
полупространству, содержащему множество $conv\ F_{1}(K)$. По условию
найдется луч $l_{1}$, такой что $\angle lOl_{1} < \varepsilon$.
Согласно ``$\varepsilon$-условию острого угла'' найдется на луче
$l_{1}$ точка $x_{1} \in l_{1}\cap K$ , такая что угол $\angle
x_{1}Oy < \frac{\pi}{2}-\varepsilon$ для всех точек $y\in
F_{1}(x_{1})$. С одной стороны, множество $F_{1}(x_{1})\subset
F_{1}(K)\subset conv F_{1}(K)$, а с другой $\angle lOy = \angle
x_{1}Oy + \angle lOl_{1} < \frac{\pi}{2}-\varepsilon + \varepsilon =
\frac{\pi}{2}$. Полученное противоречие доказывает теорему.

{\bf 5. ``$\delta$-условие коострого угла''.} Скажем, что
отображение $F$  удовлетворяет ``$\delta$-условию коострого угла'',
если для каждой точки $y^{\ast}$   некоторой $\varepsilon$-сети
$\Sigma$ на единичной сфере $S^{\ast} = \{y^{\ast} \in Y^{\ast}:
\|y^{\ast}\| = 1\}$ в $Y^{\ast}$ существует точка $x\in X$ такая,
что выполнено условие ${\rm Re}\, \langle y, y^{\ast} \rangle >
\delta \|y\|$ для всех точек $y\in F(x)$.

\newpage

{\bf Теорема 4.} {\it Пусть $D$ --- область евклидова пространства
$E^{n} = X$. Пусть $K \subset \overline{D}$ --- подмножество в
замыкании этой области и пусть существует такое сужение $F_{1}$
многозначного отображения $F: \overline{D}\rightarrow E^{n} = Y$  на
подмножество $K$, которое удовлетворяет ``$\delta$-условию коострого
угла'' для некоторой $\frac{\varepsilon}{2}$-сети $\Sigma$ в
$S^{\ast}$, такой что $\delta > sin\ \delta >
\frac{\varepsilon}{2}$. Тогда, если $conv\ F_{1}(K) \subset
F(\overline{D})$, то $0 \in F(\overline{D})$.}

\textbf{\textit {Доказательство.}} Возьмем произвольную точку
$y_{1}^{\ast}\in S^{\ast} \subset Y^{\ast}$. Исходя из условия
теоремы, существуют точки $y^{\ast}\in \Sigma\subset S^{\ast}
\subset Y^{\ast}$, $\|y^{\ast}- y_{1}^{\ast}\| <
\frac{\varepsilon}{2}$, и $x\in X$ такие, что выполнено условие
${\rm Re}\, \langle y, y^{\ast} \rangle = \|y\|\
\|y^{\ast}\|cos\angle yOy^{\ast} = \|y\|cos\angle yOy^{\ast}
> \delta \|y\| > \|y\|sin\ \delta > \frac{\varepsilon\|y\|}{2}$ для
всех точек $y\in F_{1}(x)$. Тогда ${\rm Re}\, \langle
\frac{y}{\|y\|}, y_{1}^{\ast} \rangle = {\rm Re}\, \langle
\frac{y}{\|y\|}, y^{\ast} \rangle + {\rm Re}\, \langle
\frac{y}{\|y\|},  y_{1}^{\ast} - y^{\ast} \rangle >
\frac{\varepsilon}{2} -  \|y_{1}^{\ast} - y^{\ast}\| > 0$. Теперь
данный результат следует из теоремы 2.

{\bf Замечание 1.} Все предыдущие результаты будут справедливы, если
отображение области имеет сужение, удовлетворяющее условиям теоремы.

Исследования этой работы частично поддержаны грантом Тюбитек-НАНУ
номер 110T558.

Авторы признательны профессору К.Н.Солтанову за обсуждение
результатов и ценные замечания.


\begin{thebibliography}{99}

\bibitem{} \textsc{Зелинский~Ю.\,Б., Клищук~Б.\,А., Ткачук~М.\,В.} {\sl Теоремы о неподвижной
точке для многозначных отображений} // Збірник праць Інституту
математики НАНУ. --- 2012. --- {\bf 9}, №2. --- С. 175-179.

\bibitem{}  \textsc{Красносельский~М.\,А.} {\sl Топологические методы в теории нелинейных интегральных уравнений.} --- Москва: Гостехиздат. --- 1956. --- 392 с.

\bibitem{}  \textsc{Солтанов~К.\,Н.} {\sl  О нелинейных отображениях и разрешимости нелинейных уравнений} // Докл. АН СССР. --- 1986. --- {\bf 289}, № 6. --- С. 1318---1323.

\bibitem{}  \textsc{Soltanov~K.\,N.} {\sl Remarks on Separation of Convex Sets, Fixed-Point Theorem and Applications in Theory of Linear
Operators} // Fixed Point Theory and Applications. --- 2007.
--- 14 p.

\bibitem{} \textsc{Soltanov~K.\,N.} {\sl On semi-continuous mappings, equations and inclusions in  the Banach space} //
Hacettepe J. Math. Statist. --- 2008. --- {\bf 37}. --- P.~9---24.

\bibitem{} \textsc{Зелинский~Ю.\,Б.} {\sl  Многозначные отображения в анализе.} --- Киев: Наук. думка. --- 1993. --- 264 с.

\end{thebibliography}
\end{document}